\documentclass[12pt,a4paper,twoside]{article}

\usepackage[latin1]{inputenc}
\usepackage[brazil,english]{babel}
\usepackage{amsfonts}
\usepackage{amssymb}
\usepackage{amsmath}
\usepackage{amsthm}
\usepackage[usenames,dvipsnames]{color}
\usepackage{enumitem}
\usepackage[normalem]{ulem}
\usepackage{setspace}

\numberwithin{equation}{section}
\theoremstyle{plain} 
\newtheorem{theorem}{Theorem}[section]
\newtheorem*{theorem*}{Theorem}

\newtheorem*{lemma*}{Lemma}

\theoremstyle{definition} 

\theoremstyle{remark} 
\newtheorem{remark}[theorem]{Remark}

\makeatletter
\newcommand{\aslabel}[1]{#1\def\@currentlabel{#1}}
\newcommand{\nowlabel}[1]{\def\@currentlabel{#1}}
\makeatother

\newcommand{\R}{{\mathbb{R}}}

\newcommand{\Cont}{{\mathcal{C}}}

\newcommand{\tref}[1]{\ref{#1}}
\newcommand{\pref}[1]{(\ref{#1})}

\newcommand{\pt}[1]{\left({#1}\right)}

\newcommand{\pg}[1]{\left\{{#1}\right\}}
\newcommand{\n}[1]{\left\|{#1}\right\|}

\usepackage{pstricks}\usepackage{pst-xkey}\usepackage{pst-3dplot}

\renewcommand{\l}{\ell}
\newcommand{\uu}{w}

\begin{document}
\title{On necessary conditions for the Comparison Principle and the Sub and Supersolutions Method for the stationary Kirchhoff Equation}


\author{ \sc{Leonelo Iturriaga}\,$^a$\ \  and\ \ {\sc Eugenio Massa}\,$^{b,}$\footnote{Corresponding author: tel ++55 16 33736635, fax ++55 16 33739650.}
\\[0.4cm]
$^a$ {\small Departamento de Matem\'atica, Universidad T\'ecnica Federico Santa Mar\'ia}
\\[-0.2cm] {\small Avenida Espa\~{n}a 1680, Casilla 110-V, Valpara\'{\i}so, Chile }
\\[-0.2cm] {\small email: leonelo.iturriaga@usm.cl}
\\[0.2cm]
$^b$ \small Departamento de Matem\'atica,
\\[-0.2cm]\small Instituto de Ci\^encias Matem\'aticas e de Computa\c c\~ao, Universidade de S\~ao Paulo,
\\[-0.2cm]\small Campus de S\~ao Carlos, Caixa Postal 668, 13560-970, S\~ao Carlos SP, Brazil.
\\[-0.2cm]\small e-mail: eug.massa@gmail.com 
}

\date{}

\maketitle

%

\begin{abstract}
In this paper we propose a counterexample to the validity of the Comparison Principle and of the Sub and Supersolution Method for nonlocal problems like the stationary Kirchhoff Equation. This counterexample shows that  in general smooth bounded domains  in any dimension, these properties cannot hold true if the nonlinear nonlocal term $M(\n{u}^2)$ is somewhere increasing
with respect to the $H_0^1$-norm of the solution.

Comparing with existing results, this fills a gap between known conditions on $M$ that guarantee or prevent these properties, and leads to a  condition which is necessary and sufficient  for the validity of the Comparison Principle.
%
%
%
%

\end{abstract}

\noindent  {\it  Mathematical Subject Classification MSC2010:} 
35J25   
(35B51)   

\noindent {\it Key words and phrases:} Nonlocal elliptic problems, Kirchhoff equation, Comparison Principle, Sub and Supersolution Method, counterexamples.

\section{Introduction}
In this paper we consider  problems in the form 
$$
\begin{cases}-M(\n{u}_H^2)\Delta u = f (x, u)&in\ \ \Omega\,, \\
u=0&on \ \ \partial \Omega\,,
\end{cases}\leqno{(K)}
$$
where $\Omega\subset\mathbb{R}^{N}$ is a bounded and smooth domain,  $M$ is a nonnegative function, $\n{\cdot}_H$ is the norm in $H^1_0(\Omega)$ and $f$ is some nonlinearity.

The main feature of this problem is the presence of the term  $M(\n{ u}_H^2)$, which is said to be non-local, for depending not only on the point in $\Omega$  where the equation is evaluated, but on the norm of the whole solution.

The equation in (K) is usually called of Kirchhoff type, actually a famous and important example  is the (stationary) Kirchhoff equation, originally proposed in \cite{Kirch1883} as an improvement of the vibrating string equation, in order to take into account the variation  in the  tension of the string due to the variation in  its  length with respect to the unstrained position. In the  Kirchhoff case the proposed function  $M$ takes the form  $M(t)=a+bt$ with $a,b>0$, however, different functions can be considered, either for modeling different physical phenomena (see examples in \cite{Vill_KirchhAplic}), or  for the pure mathematical interest of the problem.
For more  recent literature about  Kirchhoff type equations  like (K), we  cite  the works 
 \cite{AlvCorMa,Ma_kirchhSurv,CorFig_Kirchh,ChengWuLiu_kirch_conc,TaCh16_kirchh_GtStSiChg,SoChTa16_kirchreshigh},   which deal with the existence of  solutions with various types of nonlinearities $f$ and use mainly  variational methods.

It is also  worth noting that equations similar to (K), but with the  $p(x)$-Laplacian operator in the place of the Laplacian, have gained  interest recently for appearing in models of thermo-convective flows of non-Newtonian fluids or of electrorheological fluids, and in  image restoring problems (see in \cite{DaiHao09_p(x)Lap} and the references therein). Other nonlocal problems, though of different form with respect to (K), have great importance for their applications and also for the mathematical challenges they present: we cite for example the   ``mean field equation"
$$-\Delta u=\frac{e^u}{\int_\Omega e^u}\,,$$ which appears in the study of turbulent flows modeled by  the  Euler equation: see  \cite{CaLiMaPu_meanfield}.

\par \medskip
    
The aim of this paper is related to the recent paper   \cite{ItuGM_contra}, where the authors  obtained a nice counterexample to the validity of the Comparison Principle and of the Sub and Supersolution Method (CP and SSM for short) for the stationary Kirchhoff Equation.  This counterexample shows  how the presence of the  nonlocal term may have the effect of making unavailable these techniques, largely used in the local case.  
The motivation for  such counterexamples was to clear up some results in literature which claimed (and used) the validity of 
CP and  SSM 
for problem (K) or for  its generalizations involving the $p-$Laplacian or the $p(x)-$Laplacian, under the assumption that  $M\geq m_0>0$ is an increasing function (we refer to \cite{ItuGM_contra}  for more details and references on this matter).
In fact, the result in \cite{ItuGM_contra} shows that if the function $M$ increases fast enough, then the CP  (both in its weak and strong form) and the SSM can not hold true, at least in a ball. 

It is known that  CP and SSM hold true for the local problem where $M$ is constant, but also for the nonlocal problem (K) if $M\geq0$ is  nonincreasing but the product  $M(t^2) t$ is  increasing, as proved in \cite{AlvCor_Kirch2}.
The condition that $M(t^2) t$ is  increasing is easily seen to be necessary for the Comparison Principle to hold.
Actually, if $M(t_1^2)t_1\geq M(t_2^2)t_2$ for some positive $t_1<t_2$, taking  $\phi_1$ as the first eigenfunction of the Laplacian in a set $\Omega$, normalized with $\n{\phi_1}_H=1$, then the functions   $\l=t_2\phi_1$ and $\uu=t_1\phi_1$  satisfy 
\begin{equation*}
\begin{cases}
-M(\n{\l}_H^2)\Delta \l=M(t_2^2)t_2\lambda_1\phi_1\leq M(t_1^2)t_1\lambda_1\phi_1=  -M(\n{\uu}_H^2)\Delta \uu&in\ \ \Omega\,, \\
\l=\uu=0 &on \ \ \partial \Omega\,,
\end{cases}
\end{equation*} 
but $\l>\uu$ in $\Omega$. 

On the other hand, the result in \cite{ItuGM_contra} excludes the validity of CP and SSM if $M$ increases enough, however the authors  only consider the problem in a ball and they  do not address the question whether some growth condition on the nonlocal term could be enough to make these properties available. In particular, it is worth noting that  the hypothesis on $M$ required  for their counterexample to  work is such that the dimension $N$ has to be at least $3$, and moreover   the original Kirchhoff nonlocal term  $M(t)=a+bt$  would satisfy their condition only for $N\geq5$.

In this paper we obtain a new counterexample  that holds in any dimension, in general smooth bounded domains, and that fills the gap between the results in \cite{ItuGM_contra} and those in \cite{AlvCor_Kirch2}, showing that a necessary condition for CP and SSM to hold true in their standard form  is that the function $M$ is nonincreasing.

In fact, we prove  the following
\begin{theorem}\label{th_contraex_compar_Omegaqq}
Let $\Omega$ be a smooth bounded domain in $\R^N$. Suppose $M$ is not nonincreasing, that is, there exist positive $t_1<t_2$ such that $M(t_1)<M(t_2)$.
Then the Comparison Principle (both in its weak and strong form) and the Sub and Supersolution Method do not hold in $\Omega$, for the operator
 $$-M(\n{u}_H^2)\Delta u\,.$$ 
\end{theorem}

In particular, our result shows that even in the original Kirchhoff model, where $N=1$, $M(t)=a+bt$,
 and the operator only involves the Laplacian instead of the $p(x)$-Laplacian, 
the results in literature claiming the validity of  CP and SSM in their standard form can not hold true, and so their consequence should be questioned too. In fact, the simple assumption that $M$ is increasing somewhere is enough to make CP and SSM hold false.

Comparing with the results in \cite{AlvCor_Kirch2}, Theorem \tref{th_contraex_compar_Omegaqq} implies that the condition that  $M$ is nonincreasing is in fact necessary for both CP and SSM to hold true. Moreover, the two conditions that  $M$  is nonincreasing and $M(t^2)t$ is increasing turn out to be  necessary and sufficient, at least for the validity of the Comparison Principle.

\begin{remark}\normalfont
It is worth noting that in 
\cite{AlvCor_Kirch} a Sub and Supersolutions Method is developed, which  can deal with problem (K) and an increasing function $M$. The result is obtained by using a kind of Minty-Browder Theorem for a suitable  pseudomonotone operator, but in place of the subsolution the authors need to assume the existence of a whole family of functions which satisfy a stronger condition than just being  subsolutions: this stronger condition restricts the possible right hand sides in (K), so that, for instance, it could not be satisfied for problem \pref{eq_linThet_th} below.

Another Sub and Supersolutions Method for nonlocal problems is obtained in  \cite{AlvCov15_kirchsusup_uq} for a  problem with a nonlocal term containing a Lebesgue norm, instead of the Sobolev norm that appears in (K).
\end{remark}

\section{Proof of the result}
In order to prove that the Sub and Supersolution Method does not hold, we will provide two functions $\underline u,\overline u\in H_0^1(\Omega)\cap \Cont (\Omega)$ and a number $\Theta>0$ such that (in weak sense)
\begin{equation}\label{eq_subsup_th}\begin{cases}
 -M(\n{\underline u}_H^2)\Delta \underline u\,\leq\, \Theta \underline u&in\ \ \Omega\,, \\ 
 -M(\n{\overline u}_H^2)\Delta \overline u\,\geq\, \Theta \overline u&in\ \ \Omega\,, \\ 
 \underline u\,\leq\,\overline u&in\ \  \Omega\,,\\
 \underline u\,\leq\,0\,\leq\, \overline u&on\ \ \partial \Omega\,,\\
\end{cases}
\end{equation}
 but there exists no solution $u$ satisfying $\underline u\leq u\leq \overline u$, 
for the problem 
  \begin{equation}\label{eq_linThet_th}\begin{cases}
 -M(\n{u}_H^2)\Delta u=\Theta u&in\ \ \Omega\,, \\ 
 u=0&on \ \ \partial \Omega\,.\\
 \end{cases}
\end{equation}

The counterexample to the Weak (resp. Strong) Comparison Principle will be obtained by providing  functions $\l,\uu\in H_0^1(\Omega)\cap \Cont (\Omega)$ such that (in weak sense)
 \begin{equation}\label{eq_comp_th}
\begin{cases}
 -M(\n{\l}_H^2)\Delta \l\,<\, -M(\n{\uu}_H^2)\Delta \uu&in\ \ \Omega\,, \\ 
 \l\,\leq\, \uu&on\  \ \partial \Omega\,,\\
 \end{cases}
\end{equation}
but there exists $p\in\Omega$ where $\l(p)>\uu(p)$ (resp,  where $\l(p)=\uu(p)$).

\begin{remark}\normalfont
As one can see in the proof below, the functions $\l$ and $\underline u$ are in fact smooth, while $\uu$ and $\overline u$ are built as the gluing of two smooth functions.  By approximating  with  $\Cont^2$ functions one can modify the  counterexample so that all the functions involved are smooth and the equations are satisfied in classical sense.

It would also be possible to generalize the counterexample to the case where one considers the $p$-Laplacian instead of the Laplacian and the $W_0^{1,p}$-norm instead of the $H_0^1$-norm  in (K): actually, in the proof  we only exploit the homogeneity of the operator, not its linearity. 
\end{remark}

The idea behind the three counterexamples is similar: since the Comparison Principles hold true for the Laplacian, it is necessary that $\l,\uu$ satisfying \pref{eq_comp_th} do not satisfy $-\Delta \l\leq-\Delta \uu$. Then one has to exploit   the larger value of $M(t_2)$  with respect to $M(t_1)$ in order 
to revert the inequality and have \pref{eq_comp_th}  satisfied. Since $t_1<t_2$ this means that the $\n{\uu}_H$ must be larger than $\n{\l}_H$, which is obtained by choosing  the function $\uu$ with a large gradient near the boundary.

The same strategy is used for the counterexample to the Sub and Supersolution Method. 
In this case, inspired by the similar counterexample in  \cite{ItuGM_contra}, we  considered  a linear right hand side in \pref{eq_linThet_th} so that all possible positive solutions are known.  We  build the (strict) subsolution  and supersolution in such a way that they touch at the origin: this  only leaves one possible solution  between them. At this point,  we can choose the coefficient $\Theta$ in such a way that this is not a solution of \pref{eq_linThet_th}.

The choice of the form of the two functions $\l,\uu$ (resp. $\underline u,\overline u$) was mainly led by the need to have the first  two inequalities in \pref{eq_subsup_th} satisfied. Actually, the simplest way to obtain this is to use  first eigenfunctions of the Laplacian, in particular the function $\uu$ in  \pref{eq_comp_th} (resp. $\overline u$ in \pref{eq_subsup_th}) is defined  as the minimum between a large multiple of the first eigenfunction in $\Omega$ (which provides the required "large" H-norm) with a first eigenfunction in a slightly larger set (see equation \pref{eq_uet_min}).

Below we give the proof of our result.
\begin{proof}[Proof of theorem \tref{th_contraex_compar_Omegaqq}]
First, we define a  set slightly larger than $\Omega$ as follows:    $$\Omega^\tau=\pg{x\in\R^N:\ d(x,\Omega)<\tau}\,,$$
where $d(\cdot,\cdot)$ denotes the distance in $\R^N$ and $\tau>0$ will be taken small enough, so that $\Omega^\tau$ is still a smooth domain. 

Let  $\phi^\tau,\lambda^\tau$ be the (positive) first eigenfunction and eigenvalue for the Laplacian in $\Omega_\tau$, normalized with $\n{\phi^\tau}_\infty=1$. Let also $\phi_1,\lambda_1$ be those in $\Omega$,  again normalized with $\n{\phi_1}_\infty=1$.

Observe that  when $\tau\searrow0$ one has $\Omega^\tau\searrow\Omega$. As a consequence of the variational characterization of the first eigenvalue, this implies that  $\lambda^\tau \nearrow \lambda_1$.

Then we can fix $\tau>0$ such that  
\begin{equation}
1\geq\frac{\lambda^\tau}{\lambda_1}>\frac{M(t_1)}{M(t_2)}\,
\end{equation}
and $\Theta$   such that 
\begin{equation}\label{eq_defTh}
\lambda_1M(t_1)<\Theta<\lambda^\tau M(t_2)\,.
\end{equation}

Let now  $$c_\tau=\inf\pg{t:\ \phi_1<t\phi^\tau|_{\Omega}\ \ in\ \Omega}\,:$$
 one sees that $c_\tau\geq1$ since, by the normalization, $\phi_1\geq\phi^\tau$ at the point where $\phi_1$ attains its maximum. Also, since $\phi^\tau$ is bounded away from zero in $\Omega$, the two functions $c_\tau\phi^\tau$ and $\phi_1$ must touch somewhere in $\Omega$, that is, there exists $\widetilde p\in\Omega$ such that $c_\tau\phi^\tau(\widetilde p)=\phi_1(\widetilde p)$.

For $\varepsilon\in(0,1]$, we define the set $\Omega_{\varepsilon}^{\tau}$ and the function $u_{\varepsilon,\tau}$ as follows:
$$\Omega_{\varepsilon}^{\tau}=\pg{x\in\Omega:\ \frac1\varepsilon\phi_1(x)< c_\tau\phi^\tau(x)}\,,$$ 
\begin{eqnarray}
 u_{\varepsilon,\tau}&=&\begin{cases}
c_\tau\phi^\tau(x)& if\ x\in(\Omega_{\varepsilon}^{\tau})^C\,,\\
\frac1\varepsilon\phi_1(x)&if\ x\in\Omega_{\varepsilon}^{\tau}\,,
\end{cases}
\end{eqnarray}
or, which is the same,  
\begin{equation}\label{eq_uet_min}
u_{\varepsilon,\tau}=\min\pg{c_\tau\phi^\tau,\frac1\varepsilon\phi_1}\,.
\end{equation}

Observe that, by its definition, $u_{\varepsilon,\tau}$ is a continuous function in $\Omega$,  $u_{1,\tau}\equiv\phi_1$ and the following items  hold true:
\begin{itemize}
\item  for every  $\varepsilon\in(0,1]$, one has 
\begin{equation}\label{eq_ordemuphi}\begin{cases}
u_{\varepsilon,\tau}\geq \phi_1\,& \text{in $\Omega\,$,}
\\u_{\varepsilon,\tau}(\widetilde p)= \phi_1(\widetilde p)&\text{at the point $\widetilde p\in\Omega$ where $c_\tau\phi^\tau(\widetilde p)=\phi_1(\widetilde p)\,;$}
\end{cases}
\end{equation} 
\item 
$$\begin{cases}
-\Delta u_{\varepsilon,\tau}=\lambda^\tau c_\tau\phi^\tau\geq \lambda^\tau u_{\varepsilon,\tau}&\qquad {in\ \  {(\Omega^\tau_\varepsilon)^C}\,,}
\\
-\Delta u_{\varepsilon,\tau}=\lambda_1\frac{\phi_1}\varepsilon\geq \lambda^\tau u_{\varepsilon,\tau} &\qquad {in\ \  {\Omega^\tau_\varepsilon}\,,}
\end{cases}$$
and then 
\begin{equation}\label{eq_estSuper}
-\Delta u_{\varepsilon,\tau}\geq \lambda^\tau u_{\varepsilon,\tau}\qquad\text{$in\ \ {\Omega}$ \quad(in weak sense);}
\end{equation}
\item\begin{equation}\label{eq_Ntoinfty}
\n{u_{\varepsilon,\tau}}_H^2=\frac1{\varepsilon^2}\int_{\Omega^\tau_\varepsilon}|\nabla\phi_1|^2+c_\tau^2\int_{(\Omega^\tau_\varepsilon)^C}|\nabla\phi^\tau|^2\to+\infty \qquad \text{for $\varepsilon\to0$}\,.
\end{equation}

In order to prove \pref{eq_Ntoinfty}, let $\zeta$ and $\delta_\tau$ be such that  $|\nabla\phi_1|<\zeta$ and $\phi_\tau>\delta_\tau>0$ in $\overline \Omega$. Then  for $p\in\Omega$ we may estimate $\phi_1(p)-0<\zeta d(p,\partial\Omega)$.
This implies that $$\omega^\tau_\varepsilon:=\pg{x\in\Omega: d(x,\partial\Omega)<\varepsilon \delta_\tau/\zeta}\subseteq\Omega_{\varepsilon}^{\tau}\,.$$

Moreover, 
 we may estimate 
$$|\omega^\tau_\varepsilon|\geq|\partial\Omega|\,\frac{\varepsilon \delta_\tau}{2\zeta} \qquad \text{for $\varepsilon>0$ small enough}$$ (here $|\omega^\tau_\varepsilon|$ and $|\partial\Omega|$ denote, respectively, the $N$-dimensional and the $(N-1)$-dimensional measure of the two sets), actually, taking into account the smoothness of $\Omega$, for small $\varepsilon$, the set $\omega^\tau_\varepsilon$ is just  a  smooth band of width $\varepsilon \delta_\tau/\zeta$ near $\partial \Omega$. 

Finally,  by the properties of the first eigenfunction, there exists $a>0$ such that $|\nabla\phi_1|>a$ in $\partial\Omega$, which implies that, 
 for  $\varepsilon$ small enough, $|\nabla\phi_1|>a/2$  in  $\omega^\tau_\varepsilon$.
As a consequence, 
 $$\int_{\Omega^\tau_\varepsilon}|\nabla\phi_1|^2\geq \frac{a^2|\partial\Omega| \varepsilon\delta_\tau }{8\zeta} \qquad \text{for $\varepsilon>0$ small enough},$$ and then \pref{eq_Ntoinfty} follows.
\end{itemize}
Now let $$\l_\alpha=\alpha\phi_1\,;$$ 
then
\begin{eqnarray}
\label{eq_laplL} -\Delta \l_\alpha&=&\lambda_1\alpha\phi_1=\lambda_1\l_\alpha\,,
\\\label{eq_normL}\n{\l_\alpha}_H^2&=&\alpha^2\n{\phi_1}_H^2\,.
\end{eqnarray}

For our first counterexample  we set $\alpha=1$, and
we choose $A>0$ such that    $\n{A{\l_\alpha}}^2_H=t_1$, then we
choose $\varepsilon$ such that   $\n{Au_{\varepsilon,\tau}}_H=t_2$: this is possible by continuity, actually $u_{1,\tau}=\phi_1$, so $\n{Au_{1,\tau}}_H^2=\n{A{\l_\alpha}}^2_H=t_1<t_2$, and by \pref{eq_Ntoinfty} $\n{Au_{\varepsilon,\tau}}_H\to \infty$ when $\varepsilon\to0$.

We obtain, by \pref {eq_estSuper}, \pref{eq_defTh} and \pref{eq_laplL},
\begin{equation}\label{eq_contas_subsup}
\begin{cases}
-M(\n{Au_{\varepsilon,\tau}}_H^2)\Delta Au_{\varepsilon,\tau}\geq M(t_2)\lambda^\tau Au_{\varepsilon,\tau} >\Theta Au_{\varepsilon,\tau}\,&in\ \Omega\,,
\\-M(\n{A\l_\alpha}_H^2)\Delta A\l_\alpha=M(t_1)\lambda_1A\l_\alpha<\Theta A\l_\alpha\,&in\ \Omega\,,
\end{cases}
\end{equation}
so $\overline u:=Au_{\varepsilon,\tau}$  and $\underline u:=A\l_\alpha$ satisfy all the conditions in \pref{eq_subsup_th}, moreover $\underline u(\widetilde p)= \overline u(\widetilde p)$ 
(see \pref{eq_ordemuphi}).

We only have to prove that no  solution of  problem \pref{eq_linThet_th} exists between $\underline u$ and $\overline u$. Actually, such a solution must be positive and then it must be a multiple of $\phi_1$,   but since $\underline u(\widetilde p)= \overline u(\widetilde p)$ then the only possible  choice would be $\underline u$ itself, which is not a solution by the strict inequality in \pref{eq_contas_subsup}.

Observe that the same functions  $\uu:=Au_{\varepsilon,\tau}$ and $\l:=A\l_\alpha$ satisfy \pref{eq_comp_th} with $\l(\widetilde p)=\uu(\widetilde p)$, and then provide a counterexample to the strong form of the  Comparison Principle. 

In order to get a counterexample to the weak form of the  Comparison Principle we only have to choose $\alpha>1$, so that we will have  $\l_\alpha(\widetilde p)>u_{\varepsilon,\tau}(\widetilde p)$, 
then again choose in sequence $A>0$ such that    $\n{A{\l_\alpha}}^2_H=t_1$, then  $\varepsilon$ such that   $\n{Au_{\varepsilon,\tau}}_H=t_2$: again this is possible in view of \pref{eq_Ntoinfty} and the fact that now  $\n{Au_{1,\tau}}_H^2=\n{A{\l_\alpha}}^2_H/\alpha^2<t_1<t_2$ (see \pref{eq_normL}).
\end{proof}

\begin{remark}\normalfont
In dimension one, the construction of the functions used in the counterexamples is quite straightforward and can be done explicitely. In fact, if $\Omega=\pt{-\pi/2,\pi/2}$,  one has 
$$\l_\alpha=\alpha\cos(x)\,,$$
$$u_{\varepsilon,\tau}=\min\pg{\cos\pt{\frac x L},\frac1\varepsilon\cos(x)}\,,$$
where $L=1+\frac{2}\pi\tau$.
\end{remark}

\section*{Acknowledgement}

\noindent L. Iturriaga gratefully acknowledges financial support from Programa Basal PFB 03, CMM, U. de Chile; Fondecyt grant 1161635  and USM Grant 116.12.1.

\vspace{.3cm}

\noindent E. Massa was  supported by: grant $\#$2014/25398-0, São Paulo Research Foundation (FAPESP) and  grant $\#$308354/2014-1, CNPq/Brazil. 

\bibliographystyle{amsalpha}

  \bibliography{bibfile}


\end{document}